\documentclass{amsart}

\usepackage{amsmath,amsthm,amsfonts,amscd,amssymb}

\newtheorem{thm}{Theorem}[section]

\newtheorem{lem}[thm]{Lemma}

\theoremstyle{definition}

\theoremstyle{remark}
\newtheorem{rem}{Remark}[section]

\begin{document}

\title{Totally isotropic subspaces of small height in quadratic spaces}
\author{Wai Kiu Chan, Lenny Fukshansky, and Glenn R. Henshaw}\thanks{The second author was partially supported by NSA Young Investigator Grant \#1210223 and Simons Foundation grants \#208969, 279155.}

\address{Department of Mathematics and Computer Science, Wesleyan University, Middletown, CT 06459}
\email{wkchan@wesleyan.edu}
\address{Department of Mathematics, 850 Columbia Avenue, Claremont McKenna College, Claremont, CA 91711}
\email{lenny@cmc.edu}
\address{Department of Mathematics, Engineering and Computer Science, LaGuardia Community College, 31-10 Thomson Avenue, Long Island City, NY 11101}
\email{ghenshaw5974@gmail.com}

\subjclass[2010]{Primary 11G50, 11E12, 11H06}
\keywords{heights, quadratic forms}

\begin{abstract}
Let $K$ be a global field or~$\overline{\mathbb Q}$, $F$ a nonzero quadratic form on $K^N$, $N \geq 2$, and $V$ a subspace of $K^N$.  We prove the existence of an infinite collection of finite families of small-height maximal totally isotropic subspaces of $(V,F)$ such that each such family spans $V$ as a $K$-vector space. This result generalizes and extends a well known theorem of J. Vaaler~\cite{vaaler:smallzeros2} and further contributes to the effective study of quadratic forms via height in the general spirit of Cassels' theorem on small zeros of quadratic forms. All bounds on height are explicit.
\end{abstract}

\maketitle

\def\A{{\mathcal A}}
\def\AA{{\mathfrak A}}
\def\B{{\mathcal B}}
\def\C{{\mathcal C}}
\def\D{{\mathcal D}}
\def\E{{\mathcal E}}
\def\F{{\mathcal F}}
\def\Ff{{\mathfrak F}}
\def\G{{\mathcal G}}
\def\x{{\mathcal H}}
\def\I{{\mathcal I}}
\def\J{{\mathcal J}}
\def\K{{\mathcal K}}
\def\kk{{\mathfrak K}}
\def\L{{\mathcal L}}
\def\LL{{\mathfrak L}}
\def\M{{\mathcal M}}
\def\O{{\mathcal O}}
\def\W{{\omega}}
\def\CC{{\mathfrak C}}
\def\mm{{\mathfrak m}}
\def\MM{{\mathfrak M}}
\def\OO{{\mathfrak O}}
\def\P{{\mathcal P}}
\def\R{{\mathcal R}}
\def\s{{\mathcal S}}
\def\V{{\mathcal V}}
\def\X{{\mathcal X}}
\def\XX{{\mathfrak X}}
\def\Y{{\mathcal Y}}
\def\Z{{\mathcal Z}}
\def\H{{\mathcal H}}
\def\cee{{\mathbb C}}
\def\pee{{\mathbb P}}
\def\que{{\mathbb Q}}
\def\real{{\mathbb R}}
\def\zed{{\mathbb Z}}
\def\hyp{{\mathbb H}}
\def\aaa{{\mathbb A}}
\def\ff{{\mathbb F}}
\def\Nn{{\mathbb N}}
\def\kk{{\mathfrak K}}
\def\qbar{{\overline{\mathbb Q}}}
\def\kbar{{\overline{K}}}
\def\ybar{{\overline{Y}}}
\def\kkbar{{\overline{\mathfrak K}}}
\def\ubar{{\overline{U}}}
\def\eps{{\varepsilon}}
\def\ahat{{\hat \alpha}}
\def\bhat{{\hat \beta}}
\def\gt{{\tilde \gamma}}
\def\h{{\tfrac12}}
\def\dd{{\partial}}
\def\bfa{{\boldsymbol a}}
\def\bfb{{\boldsymbol b}}
\def\be{{\boldsymbol e}}
\def\bei{{\boldsymbol e_i}}
\def\bff{{\boldsymbol f}}
\def\bc{{\boldsymbol c}}
\def\bm{{\boldsymbol m}}
\def\bk{{\boldsymbol k}}
\def\bi{{\boldsymbol i}}
\def\bl{{\boldsymbol l}}
\def\bq{{\boldsymbol q}}
\def\bu{{\boldsymbol u}}
\def\bt{{\boldsymbol t}}
\def\bs{{\boldsymbol s}}
\def\bfu{{\boldsymbol u}}
\def\bv{{\boldsymbol v}}
\def\bw{{\boldsymbol w}}
\def\bx{{\boldsymbol x}}
\def\bX{{\boldsymbol X}}
\def\bz{{\boldsymbol z}}
\def\bwy{{\boldsymbol y}}
\def\bY{{\boldsymbol Y}}
\def\bL{{\boldsymbol L}}
\def\ba{{\boldsymbol a}}
\def\bb{{\boldsymbol\beta}}
\def\bet{{\boldsymbol\eta}}
\def\bxi{{\boldsymbol\xi}}
\def\bo{{\boldsymbol 0}}
\def\bol{{\boldkey 1}_L}
\def\ep{\varepsilon}
\def\p{\boldsymbol\varphi}
\def\q{\boldsymbol\psi}
\def\rank{\operatorname{rank}}
\def\aut{\operatorname{Aut}}
\def\lcm{\operatorname{lcm}}
\def\sgn{\operatorname{sgn}}
\def\spn{\operatorname{span}}
\def\md{\operatorname{mod}}
\def\Norm{\operatorname{Norm}}
\def\dim{\operatorname{dim}}
\def\det{\operatorname{det}}
\def\Vol{\operatorname{Vol}}
\def\rk{\operatorname{rk}}
\def\ord{\operatorname{ord}}
\def\ker{\operatorname{ker}}
\def\div{\operatorname{div}}
\def\Gal{\operatorname{Gal}}
\def\GL{\operatorname{GL}}
\def\p{\operatorname{p}}
\def\q{\operatorname{q}}
\def\t{\operatorname{t}}
\def\hs{{\hat \sigma}}
\def\chr{\operatorname{char}}

\section{Introduction and statement of results}
\label{intro}

In his celebrated 1955 paper \cite{cassels:small} J. W. S. Cassels proved that an integral quadratic form which is isotropic over~$\que$ has a non-trivial integral zero of bounded height (we detail the necessary notation of heights and quadratic forms in Section~\ref{notation} below). At the heart of Cassels' argument lies a beautiful geometric idea, based on Minkowski's Linear Forms Theorem and an appropriately chosen orthogonal reflection of the quadratic space in question. The investigation of small-height zeros of quadratic forms has since been taken up by a number of authors (see \cite{cassels_overview} for an overview) with geometric ideas continuing to play a key role. In particular, geometry of numbers techniques along with orthogonal projections were used by Schlickewei~\cite{schlickewei} to prove the existence of a small-height maximal totally isotropic subspace in an isotropic rational quadratic space, and the method was then extended by Schlickewei and  W. M. Schmidt~\cite{schmidt:schlickewei} to prove the existence of a family of such small-height subspaces which generate the entire quadratic space. The Schlickewei-Schmidt results were then extended over number fields by Vaaler~\cite{vaaler:smallzeros, vaaler:smallzeros2} by the use of adelic geometry of numbers and local projection operators at all places. Vaaler's method~\cite{vaaler:smallzeros} has also been extended over global function fields of odd characteristic in~\cite{quad_zero} to prove the existence of a small-height maximal totally isotropic subspace in an isotropic quadratic space. An important underlying idea used in all such arguments over global fields is Nortcott's principle, guaranteeing finiteness of sets of projective points of bounded height and degree. Northcott's principle no longer holds over~$\qbar$, where an analogue of Vaaler's first result \cite{vaaler:smallzeros} was established in~\cite{quad_qbar} by an application of another geometric principle, namely a version of arithmetic Bezout's theorem due to Bost, Gillet, and Soul\'e~\cite{bgs}.

Estimates on small-height zeros were further applied to produce effective versions of the Witt decomposition theorem over global fields and~$\qbar$ in~\cite{witt}, \cite{quad_qbar}, and~\cite{quad_zero}. This theorem, reviewed in Section~\ref{gen_vaaler} below, states that a quadratic space can be decomposed into an orthogonal sum of its radical, a maximal anisotropic subspace, and a collection of hyperbolic planes, where all the components in the decomposition have small height.  The method employed in all coefficient field settings is the same, and  it uses the existence of a small-height maximal totally isotropic subspace as a starting point. In other words, different geometric techniques may be required to obtain estimates on small-height zeros, but once those are in place, effective structure theorems can be derived by a unified approach over different types of coefficient fields.

The main result of~\cite{vaaler:smallzeros} guarantees the existence of a small-height maximal isotropic subspaces of a quadratic space over $K$, in case $K$ is a number field, and analogous statements for global function fields and~$\qbar$ follow from~\cite{quad_zero} and~\cite{quad_qbar}, respectively. Furthermore, Vaaler's second result \cite{vaaler:smallzeros2} implies that there exists a finite family of such small-height isotropic subspaces with the additional independence condition that the corresponding subspaces span~$V$.

\begin{rem} \label{fano} Another way to view the Schlickewei-Schmidt-Vaaler type results on small-height totally isotropic subspaces is in terms of Fano varieties. The {\it Fano variety} of $m$-planes on a projective variety $X$ defined over a field $K$ is the set of $(m+1)$-dimensional vector spaces over $K$ that are contained in that projective variety; this is a natural generalization of the Grassmannian. Let $V$ be a subspace of $K^N$, $N \geq 2$, $X$ a quadratic hypersurface defined on $V$ over a field~$K$, and let $m+1$ be the dimension of a maximal linear subspace of $V$ contained in~$X$. Then the main result of~\cite{vaaler:smallzeros} guarantees the existence of a small-height point on the Fano variety of $m$-planes of~$X$ in case $K$ is a number field, and analogous statements for global function fields and~$\qbar$ follow from~\cite{quad_zero} and~\cite{quad_qbar}, respectively. Furthermore, Vaaler's second result \cite{vaaler:smallzeros2} implies that there exists a family of such small-height points in $X$ with the additional independence condition that the corresponding subspaces span~$V$.
\end{rem}

In this note we prove that in fact there exist {\em infinitely many} such families. To rigorously state our main result, we need to introduce some notation. Let $K$ be a global field or $\qbar$.  By a global field we always mean either a number field or a function field of transcendence degree 1 over a finite field $\ff_q$ of $q$ elements, where $q$ is odd. Let $V \subseteq K^N$ be an $L$-dimensional vector space, $1 \leq L \leq N$.   Suppose that $F$ is a quadratic form in $N$ variables defined over $K$.  Let $\W$ be the Witt index of the quadratic space $(V,F)$, $\lambda$ be the dimension of its radical $V^{\perp}$, and $r := L-\lambda$ be the rank of $F$ on $V$.  In particular, the dimension of a maximal totally isotropic subspace of $V$ is $\W + \lambda$. Define
\begin{equation}
\label{alpha1}
\alpha(\W,r) = \left\{ \begin{array}{ll}
\frac{(2\W + r)(\W + 1)(\W + 2)(\W + 14) + 4(\W + r + 16)}{8} & \!\!\!\!\!\mbox{if $K = $ global field},\\
(\W^2 + 1)(\W + 10)(\W+1)(\W+2)\left(\frac{3}{2}\right)^\W + \frac{\W + 14}{2} + r \ & \!\!\!\!\!\mbox{if $K = \qbar$},
\end{array} \right.
\end{equation}
and
\begin{equation}
\label{beta1}
\beta(\W) = \left\{\begin{array}{ll}
\frac{(\W + 1)(\W + 2)(\W + 14)}{2} + 1 & \mbox{ if $K = $ global field},\\
\frac{(6\W+5)(\W + 5)(\W+1)(\W+2)\left(\frac{3}{2}\right)^\W}{4\W+2} + 5 & \mbox{ if $K = \qbar$.}
\end{array} \right.
\end{equation}
For a positive integer $n$, let
\begin{equation}
\label{akn}
a_K(n) =  \left\{ \begin{array}{ll}
n^2 & \mbox{if $K = $  number field,} \\
e^{2n} & \mbox{if $K = $ function field,} \\
1 & \mbox{if $K=\qbar$.}
\end{array}
\right.
\end{equation}
In what follows below, $H$ and $\H$ are height functions defined on projective spaces.  Their definitions are reviewed in Section \ref{notation}.

\begin{thm} \label{general_vaaler}
Let $V$ be a nonzero subspace of $K^N$ and let $F$ be a quadratic form on $K^N$ of rank $r \geq 3$ on~$V$. With notation as above, there exists an infinite family $\{W_k^n\}$ of maximal totally isotropic subspaces of $V$, where $n \in \Nn$ and $k$ belongs to a finite indexing set $\I$, such that for all $n,n'  \geq 1$ and $k,k' \in \I$,
\begin{enumerate}
\item[(i)] $\W + \lambda - 4 \leq \dim_K W^n_k \cap W^{n'}_{k'} \leq \W + \lambda - 1$ if $(n,k) \neq (n',k')$;

\item[(ii)] $\{W^n_{k}: k \in \I\}$ spans $V$ over $K$ for each $n$;

\item[(iii)] $\H(W_{k}^n) \leq \CC_K(N,\W,L,r)\ a_K(n)\, H(F)^{\alpha(\W,r)}\, \H(V)^{\beta(\W)}$, for each $n$ and $k$,
\end{enumerate}
where $\CC_K(N,\W,L,r)$ is an explicit constant.
\end{thm}

The constant $\CC_K(N,\W,L,r)$ is not of great importance to the statement of the theorem, however we still derive it below in the proof of Theorem~\ref{general_vaaler} (equations \eqref{W_main_constant_nf}, \eqref{W_main_constant_qb} for $J > 0$, and equations \eqref{W0_ht_F_const_nf}, \eqref{W0_ht_F_const_qb} for $J = 0$) just to demonstrate the explicit nature of our inequalities.

Our Theorem~\ref{general_vaaler} is a generalization of Vaaler's Theorem~1 of~\cite{vaaler:smallzeros2} in the following sense. While our bound is weaker than Vaaler's, our result holds in more generality. Firstly, we are working over a much wider choice of coefficient fields. Further, we prove the existence of an {\it infinite collection} of small-height generating families of maximal totally isotropic subspaces, and the height bound depends on the index $n$ of each such family very mildly, as indicated by the value of~$a_K(n)$ (due to Northcott's finiteness property, dependence on $n$ in case of global fields is inevitable). Finally, the main exponents $\alpha(\W,r)$ and $\beta(\W)$ in our bound are independent of $n$ and are still polynomials in $r$ and $w$ in the case of global fields.

The proof of our result is different from Vaaler's argument. While the underlying starting point for both methods is the original Schlickewei-Vaaler theorem on the existence of one small-height maximal totally isotropic subspace of a quadratic space, we employ a geometric technique based on the effective Witt decomposition theorem. Our construction uses this orthogonal decomposition, along with Siegel's lemma applied to the anisotropic component, to build sequences of points spanning desired generating families of isotropic subspaces. The heights can be controlled throughout this construction, which allows for explicit estimates.

This paper is organized as follows. In Section~\ref{notation} we set the notation, define the necessary constants and height functions, and review the basic terminology in the algebraic theory of quadratic forms. We then present the proof of Theorem~\ref{general_vaaler} in Section~\ref{gen_vaaler}. Our argument works simultaneously over all coefficient fields. The main tool we use in the proof is an effective version of the Witt decomposition theorem, which we also review in Section~\ref{gen_vaaler}.

\section{Preliminaries} \label{notation}

\subsection{Notations and heights}

We start with some notation, following \cite{witt} and \cite{null}.  For the rest of this subsection,  $K$ is always a global field and $\kbar$ is its algebraic closure.    In the number field case, we write $d = [K:\que]$ for the global degree of $K$ over $\que$, $\D_K$ its discriminant, $\omega_K$ the number of roots of unity in $K$, $r_1$ its number of real embeddings, and $r_2$ its number of conjugate pairs of complex embeddings; so $d=r_1+2r_2$.

When $K$ is a function field, we will fix a $t \in K$ such that $K$ is a finite separable extension of the rational function field $\ff_q(t)$.  Its global degree is $d = [K:\ff_q(t)]$, and the effective degree of $K$ over $\ff_q(t)$ is
$$\mm(K) = \frac{[K:\ff_q(t)]}{[k_0:\ff_q]},$$
where $k_0$ is the algebraic closure of $\ff_q$ in $K$.

Let $M(K)$ be the set of all places of $K$. For each place $v \in M(K)$, we write $K_v$ for the completion of $K$ at $v$, and let $d_v$ be the local degree of $K$ at $v$, which is $[K_v:\que_v]$ in the number field case, and $[K_v:\ff_q(t)_v]$ in the function field case.   

If $K$ is a number field, then for each place $v \in M(K)$ we define the absolute value $|\ |_v$ to be the unique absolute value on $K_v$ that extends either the usual absolute value on $\real$ or $\cee$ if $v | \infty$, or the usual $p$-adic absolute value on $\que_p$ if $v|p$, where $p$ is a prime.

If $K$ is a function field, then all absolute values on $K$ are non-archimedean. For each $v \in M(K)$, let $\OO_v$ be the valuation ring of $v$ in $K_v$ and $\MM_v$ the unique maximal ideal in $\OO_v$. We choose the unique corresponding absolute value $|\ |_v$ such that:
\begin{trivlist}
\item (i) if $1/t \in \MM_v$, then $|t|_v = e$,
\item (ii) if an irreducible polynomial $p(t) \in \MM_v$, then $|p(t)|_v = e^{-\deg(p)}$.
\end{trivlist}

\noindent
In both cases, the following product formula is satisfied:
\begin{equation}
\label{product_formula}
\prod_{v \in M(K)} \vert a\vert^{d_v}_v = 1, \quad \mbox{ for all $a \in K^\times$}.
\end{equation}

If $K$ is a function field, let $g(K)$ be the genus of $K$, and for any integer $\ell \geq 1$ we define
\begin{equation}
\label{EKL}
\E_K(\ell) = e^{\frac{g(K) \ell}{d}}.
\end{equation}
Let $n(K)$ be the number of places of degree one in $K$, and $h_K$ for the number of divisor classes of degree zero of $K$.  For any integer $j \geq 1$, let
\begin{equation}
\label{RK}
R_K(j) = \frac{n(K)-1}{2} \left( (j-q+2) h_K \sqrt{n(K)} \right)^{\frac{1}{n(K)-1}} + h_K (n(K)-1) \sqrt{n(K)}.
\end{equation}

We can now define two quantities $C_K(\ell)$ and $A_K(j)$ which will appear in our various height estimates in the subsequent discussion: for any integer $\ell \geq 1$, let
\begin{equation}
\label{CKL}
C_K(\ell) = \left\{ \begin{array}{ll}
\left( \left( \frac{2}{\pi} \right)^{r_2} \vert\D_K\vert \right)^{\frac{\ell}{2d}} & \mbox{if $K = $ number field,} \\
\exp\left(\frac{(g(K)-1+\mm(K))\ell}{\mm(K)}\right) & \mbox{if $K = $ function field,} \\
e^{\frac{\ell(\ell-1)}{4}} + \eps & \mbox{if $K = \qbar$; here we can take any $\eps >0$,}
\end{array}
\right.
\end{equation}
and for any integer $j \geq 1$, let
\begin{equation}
\label{AKLM}
A_K(j) = \left\{ \begin{array}{ll}
\left( j \sqrt{2^{r_1} |\D_K|} \right)^{\frac{1}{d}} & \mbox{if $K$ is a number field with $\omega_K \leq j$,} \\
e^{R_K(j)}  & \mbox{if $K$ is a function field with $q \leq j$,} \\
1 & \mbox{otherwise.}
\end{array}
\right.
\end{equation}

If $K$ is a number field and $v \mid \infty$, then for each positive integer $j$ we define, as in \cite{vaaler:smallzeros},
\[ r_v(j) = \left\{ \begin{array}{ll}
    \pi^{-\frac{1}{2}}\, \Gamma(\frac{j}{2}+1)^{\frac{1}{j}} & \mbox{if $v \mid \infty$ is real,} \\
    (2\pi)^{-\frac12}\,\Gamma(j+1)^{\frac{1}{2j}} & \mbox{if  $v \mid \infty$ is complex.}
\end{array}
\right. \]
and
\begin{equation}
\label{constant:B}
B_K(j) = 2\vert \D_K\vert^{\frac{1}{2d}} \prod_{v \mid \infty} r_v(j)^{\frac{d_v}{d}}.
\end{equation}

For each $v \in M(K)$ we define a local height $H_v$ on $K_v^N$ by
$$H_v(\bx) = \max_{1 \leq i \leq N} \vert x_i\vert^{d_v}_v,$$
for each $\bx \in K_v^N$. Also, for each $v \mid \infty$ we define another local height
$$\H_v(\bx) = \left( \sum_{i=1}^N \vert x_i\vert_v^2 \right)^{\frac{d_v}{2}}.$$
As a result, we have two slightly different global height functions on $K^N$:
\begin{equation}
\label{global_heights}
H(\bx) = \left( \prod_{v \in M(K)} H_v(\bx) \right)^{\frac{1}{d}},\ \ \H(\bx) = \left( \prod_{v \nmid \infty} H_v(\bx) \times \prod_{v \mid \infty} \H_v(\bx) \right)^{\frac{1}{d}},
\end{equation}
for each $\bx \in K^N$. These height functions are {\it homogeneous}, in the sense that they are defined on the projective space over $K^N$ thanks to the product formula (\ref{product_formula}).  It is easy to see that
\begin{equation}
\label{ht_ineq_sqrt}
H(\bx) \leq \H(\bx) \leq \sqrt{N} H(\bx).
\end{equation}
Notice that in case $K$ is a function field, $M(K)$ contains no archimedean places, and so $H(\bx) = \H(\bx)$ for all $\bx \in K^N$. We also define the {\it inhomogeneous} height
$$h(\bx) = H(1,\bx),$$
which generalizes the Weil height on algebraic numbers.  Clearly, $h(\bx) \geq H(\bx)$ for each $\bx \in K^N$.  We extend the height functions $H$ and $h$ to polynomials by evaluating the height of its coefficient vectors, and to matrices by viewing them as vectors.  However, if $X$ is a matrix with $\bx_1, \ldots, \bx_L$ as its columns, then $\H(X)$ is always the height $\H(\bx_1\wedge \dots \wedge \bx_L)$, where $\bx_1 \wedge \dots \wedge \bx_L$ is viewed as a vector in $K^{\binom{N}{L}}$ under the standard embedding.

Let $V$ be an $L$-dimensional subspace of $K^N$.  Then there exist $N \times L$ matrix $X$ and $(N-L) \times N$ matrix $A$, both are over $K$, such that
$$V = \{ X \bt : \bt \in K^L \} = \{ \bx \in K^N : A \bx = 0 \}.$$
The Brill-Gordan duality principle \cite{gordan:1} (also see Theorem 1 on p. 294 of \cite{hodge:pedoe}) implies that $\H(X)=\H(A^t)$, and $\H(V)$ is defined to be this common value.  This coincides with the choice of heights in \cite{vaaler:siegel}.

An important observation is that due to the normalizing exponent $1/d$ in (\ref{global_heights}) all our heights are {\it absolute}, meaning that they do not depend on the number field or function field of definition.  In particular, we can extend their definitions to ~$\kbar$.

Let
\begin{equation}
\label{delta}
\delta = \left\{ \begin{array}{ll}
1 & \mbox{if char$K = 0$}, \\
0 & \mbox{if char$K > 0$}.
\end{array}
\right.
\end{equation}
We will also need a few technical lemmas detailing some basic properties of heights. The first one bounds the height of a linear combination of vectors.

\begin{lem} \label{sum_height} For $\xi_1,...,\xi_L \in K$ and $\bx_1,...,\bx_L \in K^N$,
$$H \left( \sum_{i=1}^L \xi_i \bx_i \right) \leq h \left( \sum_{i=1}^L \xi_i \bx_i \right) \leq L^{\delta} h(\bxi) \prod_{i=1}^L h(\bx_i),$$
where $\bxi = (\xi_1,...,\xi_L) \in K^L$, and $\delta$ is as in (\ref{delta}) above.
\end{lem}

\proof
The first inequality is clear from the definition of height functions. To prove the second inequality, let $v \in M(K)$. If $v \mid \infty$, then
\begin{eqnarray*}
\max \left\{ 1, H_v \left( \sum_{i=1}^L \xi_i \bx_i \right) \right\} & \leq & L \max_{1 \leq i \leq L} \left\{ \max\{ 1,|\xi_i|_v \} \max \{1, H_v(\bx_i) \} \right\} \\
& \leq & L \prod_{i=1}^L\left\{ \max\{ 1,|\xi_i|_v \} \max \{1, H_v(\bx_i) \} \right\},
\end{eqnarray*}
and if $v \nmid \infty$, then
\begin{eqnarray*}
\max \left\{ 1, H_v \left( \sum_{i=1}^L \xi_i \bx_i \right) \right\} & \leq & \max_{1 \leq i \leq L} \left\{ \max\{ 1,|\xi_i|_v \} \max \{1, H_v(\bx_i) \} \right\} \\
& \leq & \prod_{i=1}^L\left\{ \max\{ 1,|\xi_i|_v \} \max \{1, H_v(\bx_i) \} \right\}.
\end{eqnarray*}
The conclusion of the lemma now follows by taking a product over all places in $M(K)$ while keeping in mind that all absolute values over a function field are non-archimedean.
\endproof

The following two are adaptations of Lemma 4.7 of \cite{absolute:siegel} and Lemma 2.3 of \cite{witt}, respectively, to our choice of height functions, using~\eqref{ht_ineq_sqrt}.

\begin{lem} \label{lem_4.7} Let $V$ be a subspace of $K^N$, $N \geq 2$, and let subspaces $U_1,\dots,U_n \subseteq V$ and vectors $\bx_1,\dots,\bx_m \in V$ be such that
$$V = \spn_K \{ U_1,\dots,U_n,\bx_1,\dots,\bx_m \}.$$
Then
$$\H(V) \leq N^{\frac{\delta m}{2}} \H(U_1) \dots \H(U_n) H(\bx_1) \dots H(\bx_m),$$
where $\delta$ is as in (\ref{delta}) above.
\end{lem}

\begin{lem} \label{lem_2.3} Let $X$ be a $J \times N$ matrix over $K$ with row vectors $\bx_1,...,\bx_J$, and let $F$ be a symmetric bilinear form in $N$ variables over $K$ (we also write $F$ for its $N \times N$ coefficient matrix). Then
$$\H(X F) \leq N^{\frac{3J \delta}{2}} H(F)^J \prod_{i=1}^J H(\bx_i),$$
where $\delta$ is as in (\ref{delta}) above.
\end{lem}

The next one is Lemma 2.2 of \cite{witt} over any global field.

\begin{lem} \label{intersection} Let $U_1$ and $U_2$ be subspaces of $K^N$. Then
$$\H(U_1 \cap U_2) \leq \H(U_1) \H(U_2).$$
\end{lem}

\begin{rem} \label{overkbar} Lemmas~\ref{sum_height} - \ref{intersection} also hold verbatim with $K$ replaced by $\kbar$.
\end{rem}

\subsection{Quadratic Forms}

Here we introduce some basic language of quadratic forms which are necessary for subsequent discussion.  For an introduction to the subject, the readers are referred to, for instance, Chapter 1 of \cite{scharlau}.  For the sake of more generality, we allow $K$ to be any field of characteristic not 2.  We write
$$F(\bX,\bY) = \sum_{i=1}^N \sum_{j=1}^N f_{ij} X_i Y_j$$
for a symmetric bilinear form in $2N$ variables with coefficients $f_{ij} = f_{ji}$ in $K$, and $F(\bX) = F(\bX,\bX)$ for the associated quadratic form in $N$ variables; we also use $F$ to denote the symmetric $N \times N$ coefficient matrix $(f_{ij})_{1 \leq i,j \leq N}$.    Let $V$ be an $L$-dimensional subspace of $K^N$, $2 \leq L \leq N$.  Then $F$ is also defined on $V$, and we write $(V,F)$ for the corresponding quadratic space.

For each subspace $U$ of $(V,F)$, its radical is the subspace
$$U^{\perp} := \{ \bx \in U : F(\bx, \bwy) = 0\ \forall\ \bwy \in U \}.$$
We define $\lambda(U):=\dim_K U^{\perp}$, and will write $\lambda$ to denote $\lambda(V)$. A subspace $U$ of $(V,F)$ is called regular if $\lambda(U)=0$.

A point $\bo \neq \bx \in V$ is called isotropic if $F(\bx)=0$ and anisotropic otherwise. A subspace $U$ of $V$ is called isotropic if it contains an isotropic point, and it is called anisotropic otherwise. A totally isotropic subspace $W$ of $(V,F)$ is a subspace such that for all $\bx,\bwy \in W$, $F(\bx,\bwy)=0$. All maximal totally isotropic subspaces of $(V,F)$ contain $V^{\perp}$ and have the same dimension. Given any maximal totally isotropic subspace $W$ of $V$, let
$$\W = \W(V) := \dim_K(W)-\lambda,$$
which is the Witt index of $(V,F)$.  If $K=\kbar$, then $\W= [(L-\lambda)/2]$, where $[\ ]$ stands for the integer part function.

Two subspaces $U_1$ and $U_2$ of $(V,F)$ are said to be orthogonal if $F(\bx,\bwy) = 0$ for all $\bx \in U_1$ and $\bwy \in U_2$, in which case we write $U_1 \perp U_2$ for their orthogonal direct sum.  Two vectors $\bx,\bwy \in V$ are called a hyperbolic pair if $F(\bx) = F(\bwy) = 0$ and $F(\bx,\bwy) \neq 0$; the subspace $\hyp(\bx,\bwy) := \spn_K \{\bx,\bwy\}$ that they generate is regular and is called a hyperbolic plane.  It is well known that there exists a Witt decomposition of the quadratic space $(V,F)$ of the form
\begin{equation}
\label{decompose}
V = V^{\perp} \perp \hyp_1 \perp\ \dots \perp \hyp_{\W} \perp U,
\end{equation}
where $\hyp_1, \dots, \hyp_{\W}$ are hyperbolic planes and $U$ is an anisotropic subspace, which is determined uniquely up to isometry. The rank of $F$ on $V$ is $r:=L-\lambda$. In case $K=\kbar$, $\dim_K U = 1$ if $r$ is odd and 0 if $r$ is even.

\section{Families of generating isotropic subspaces}
\label{gen_vaaler}

Let $K$ be a global field or~$\qbar$, as specified in Section~\ref{intro}.  Let $F$ be a nonzero quadratic form in $N$ variables, and $V \subseteq K^N$ be an
$L$-dimensional quadratic space such that $(V,F)$ has  rank $r$ and Witt index $\W \geq 1$.   For the proof of Theorem \ref{general_vaaler}, we need a Witt decomposition (\ref{decompose}) of $(V,F)$ in which the height of every component is explicitly bounded above.    We first introduce some constants that will appear in those upper bounds.  First of all, when $K$ is a global field, let
\begin{equation}
\label{constant}
\G_K(N,L,\W) = \left\{ \begin{array}{ll}
\left\{ \left( 2^{2\W+1} B_K(L)^2 \right)^L \left( N \vert \D_K\vert^{\frac{1}{d}} \right)^{\W+5L} \right\}^{\frac{\W(\W+3)}{8}} & \mbox{if char$K =  0$}, \\
\left\{ C_K(L)^2 q^{\frac{(L^2 - \W + \W^2)g(K)}{d}} \right\}^{\frac{\W(\W+1)}{4}} & \mbox{if char$K\neq  0$},
\end{array}
\right.
\end{equation}
where $C_K(L)$ and $B_K(L)$ are defined in \eqref{CKL} and \eqref{constant:B} respectively.  Second, we define $\eta(L,r)$ by
$$\eta(L,r) = \left\{ \begin{array}{ll}
3^{\frac{L(L-1)}{2}} & \mbox{if $r < L$,} \\
1 & \mbox{if $r=L$.}
\end{array}
\right.$$

\begin{thm} \label{witt_decomp}
There exists a Witt decomposition
\begin{equation}
\label{witt_dec}
V = V^{\perp} \perp \hyp_1 \perp\ \dots \perp \hyp_{\W} \perp U
\end{equation}
for the quadratic space $(V,F)$ such that
\begin{equation}
\label{sing_height}
\H(V^{\perp}) \leq \left\{ \begin{array}{ll}
B_K(r)^r\, H(F)^{\frac{r}{2}} \,\H(V) & \mbox{if $K= $  number field,} \\
q^{\frac{rg(K)}{d}}\, H(F)^{\frac{r}{2}} \,H(V) & \mbox{if $K = $ function field,} \\
3^{\frac{L(L-1)}{2}}\, H(F)^r\, \H(V)^2  & \mbox{if $K=\qbar$.}
\end{array}
\right.
\end{equation}
Moreover, if $K$ is a global field,
\begin{equation}
\label{height}
\max\{ \H(\hyp_i),\H(U) \} \leq \G_K(N,L,\W) \left\{ H(F)^{\frac{2\W + r}{4}} \H(V) \right\}^{\frac{(\W+1)(\W+2)}{2}},
\end{equation}
for each $1 \leq i \leq \W$; and if $K = \qbar$, then
\begin{equation}
\label{witt1}
\H(\hyp_i) \leq 3^{12 \W^4 (\W+1) \left(\frac{3}{2}\right)^\W} \left\{ \sqrt{\W}\ H(F)^{\W^2+1} \left(  \eta(L,r) \H(V) \right)^{\frac{6\W+5}{4\W+2}} \right\}^{\frac{(\W+1)(\W+2)}{2} \left(\frac{3}{2}\right)^\W},
\end{equation}
and
\begin{equation}
\label{witt2}
\H(U)  \leq 2 \sqrt{2\W+1}\ 3^{\frac{(2\W+3)\W}{2}} \left(  \eta(L,r) \H(V) \right)^{\frac{2\W+3}{4\W+2}}.
\end{equation}
\end{thm}

\proof
The asserted existence of an effective Witt decomposition for a quadratic space has been established in~\cite{witt} over number fields, in~\cite{quad_qbar} over~$\qbar$, and follows from the results of~\cite{quad_zero} over global function fields of odd characteristic.

More specifically, if $K$ is a number field the statement of the theorem is precisely Theorem~1.3 of~\cite{witt}. The proof of this theorem relies on Theorem~1 of~\cite{vaaler:smallzeros}, establishing the existence of a small-height totally isotropic subspace of a regular quadratic space, and Theorem~2 of~\cite{vaaler:smallzeros2}, giving a bound on the height of the radical of a quadratic space. Analogues of these results over a global function field of odd characteristic are given by Theorem~3.1 and Lemma~3.10 of~\cite{quad_zero}, respectively. With the use of these results, the proof of Theorem~1.3 of~\cite{witt} now carries over to the case of a global function field of odd characteristic word for word, establishing the theorem in this case.

Over $\qbar$, inequalities~\eqref{witt1} and~\eqref{witt2} for a regular quadratic space are established by Theorem~5.1 of~\cite{quad_qbar}. If the space $(V,F)$ is not regular, then it can be decomposed as
\begin{equation}
\label{reg_rad}
V = V^{\perp} \perp R,
\end{equation}
where $(R,F)$ is regular and heights of $V^{\perp}$ and $R$ are bounded as in Lemma~3.5 of~\cite{quad_qbar}. Now applying Theorem~5.1 of~\cite{quad_qbar} to $(W,F)$ yields the result in the~$\qbar$ case.
\endproof

\begin{rem} \label{reg_rem} As remarked in~\eqref{reg_rad} above, the Witt decomposition~\eqref{witt_dec} can be viewed as $V=V^{\perp} \perp R$, where $R = \hyp_1 \perp\ \dots \perp \hyp_{\W} \perp U$ is $r$-dimensional maximal regular subspace of $V$ with respect to $F$. Now, the dependence on $r$ instead of $L$ in the exponent of $H(F)$ in our bound~\eqref{height} is obtained by applying the effective Witt decomposition theorem of~\cite{witt} and~\cite{quad_zero} to the regular space $(R,F)$ instead of $(V,F)$: this is possible, since $\H(R) \ll \H(V)$ (see, for instance Lemma~3.2 of~\cite{witt}).
\end{rem}

We will also need a lemma on the existence of a small-height hyperbolic pair in a given hyperbolic plane, which is Lemma~4.3 of~\cite{quad_zero} when $K$ is a global field.  The proof for the $\qbar$ case is exactly the same.

\begin{lem} \label{hyper} Let $F$ be a symmetric bilinear form in $2N$ variables over $K$. Let $\hyp \subseteq K^N$ be a hyperbolic plane with respect to $F$. Then there exists a basis $\bx,\bwy$ for $\hyp$ such that
$$F(\bx) = F(\bwy) = 0,\ F(\bx,\bwy) \neq 0,$$
and
\begin{equation}
\label{x_bound}
H(\bx) \leq h(\bx) \leq \left\{ \begin{array}{ll}
2 \sqrt{2}\ B_K(1)^2\, H(F)^{\frac{1}{2}}\, \H(\hyp) & \mbox{if $K = $ number field,} \\
q^{\frac{4g(K)}{d}}\, H(F)^{\frac12} \,\H(\hyp), & \mbox{if $K= $ function field,} \\
72\ H(F)^{\frac{1}{2}} \,\H(\hyp)^2 & \mbox{if $K = \qbar$,}
\end{array}
\right.
\end{equation}
as well as
\begin{equation}
\label{y_bound}
H(\bwy) \leq h(\bwy) \leq \left\{ \begin{array}{ll}
24 \sqrt{2} N^2\ \left( B_K(1) G_K \right)^2 H(F)^{\frac{3}{2}}\, \H(\hyp)^3 & \mbox{if $K = $  number field,} \\
4q^{\frac{4g(K)}{d}} \,G_K^2\, H(F)^{\frac{3}{2}}\, \H(\hyp)^3 & \mbox{if $K = $ function field,} \\
864\, N^2\ G_K^2\, H(F)^{\frac{3}{2}}\, \H(\hyp)^4 & \mbox{if $K = \qbar$,}
\end{array}
\right.
\end{equation}
where $G_K = \E_K(2)^{1-\delta} A_K(2) C_K(2)$ and $\delta$ as in \eqref{delta}.
\end{lem}

Let
$$V = V^\perp  \perp \hyp_1 \perp\ \dots \perp \hyp_{\W} \perp U$$
be the Witt decomposition of the $L$-dimensional quadratic space $(V,F)$ obtained from Theorem \ref{witt_decomp}. For $1\leq i \leq \W$, let $\bx_i,\bwy_i$ be a small-height hyperbolic pair for the corresponding hyperbolic plane $\hyp_i$, as guaranteed by Lemma~\ref{hyper}.   We may assume that $h(\bx_i) \leq h(\bwy_i)$, $H(\bx_i) \leq H(\bwy_i)$ for each $1 \leq i \leq \W$.  Define $J$ to be the dimension of the anisotropic component $U$, which is equal to $r - 2\W$.

If $J > 0$, let $\bu_1,\dots,\bu_J$ be the small-height basis for $U$, guaranteed by Siegel's lemma (see \cite{vaaler:siegel} and \cite{absolute:siegel}, conveniently formulated in Theorem~1.1 of \cite{null}, as well as Theorem~1.2 of \cite{null}):
\begin{equation}
\label{siegel_for_U}
\prod_{k=1}^J H(\bu_k) \leq \prod_{k=1}^J h(\bu_k) \leq C_K(J) \E_K(J)^{1-\delta} \H(U),
\end{equation}
where $\delta$ is as in \eqref{delta}.

The indexing set $\I$ that we use in the statement of Theorem~\ref{general_vaaler} is more conveniently written as the set of pairs of indices, defined by
$$\I = \left\{ \begin{array}{ll}
\{ (i,j) : 1 \leq i \leq \W,\ 1 \leq j \leq J \} & \mbox{if $J > 0$,} \\
\{ (i,j) : 1 \leq i \neq j \leq \W \}  & \mbox{if $J=0$.}
\end{array}
\right.$$
Here we are assuming that $\W>1$, and will separately treat the case $\W=1$ at the end of this section. For each pair $(i,j) \in \I$, define
\begin{equation}
\label{aij}
\alpha_{ij} = \left\{ \begin{array}{ll}
- \displaystyle{\frac{F(\bu_j)}{2F(\bx_i,\bwy_i)}} & \mbox{if $J > 0$}, \\
- \displaystyle{\frac{F(\bx_j,\bwy_j)}{F(\bx_i,\bwy_i)}}  & \mbox{if $J=0$,}
\end{array}
\right.
\end{equation}
so that $\alpha_{ij} \neq 0$. For each integer $n \geq 1$, let
\begin{equation}
\label{xin}
\xi_n = \left\{ \begin{array}{ll}
n & \mbox{if $K = $  number field,} \\
t^n & \mbox{if $K = $  function field,} \\
e^{\frac{2\pi i}{n}}  & \mbox{if $K=\qbar$, where $i=\sqrt{-1}$.}
\end{array}
\right.
\end{equation}
Now, for each pair $(i,j) \in \I$ and each $n \geq 1$, define subspace $W^n_{ij}$ of $V$ by
\begin{equation}
\label{wnij_1}
W^n_{ij} = \spn_K \left\{ V^{\perp}, \bx_1,\dots,\bx_{i-1},\bx_{i+1},\dots,\bx_{\W},\bx_i+\xi_n^2 \alpha_{ij} \bwy_i + \xi_n \bu_j \right\}
\end{equation}
when $J > 0$, and
\begin{eqnarray}
\label{wnij_2}
W^n_{ij} & = & \spn_K \Big\{ V^{\perp}, \bx_1,\dots,\bx_{i-1},\bx_{i+1},\dots,\bx_{j-1}, \nonumber \\
& &\ \ \ \ \ \ \ \ \ \ \ \bx_{j+1},\dots,\bx_{\W},\bx_i + \xi_n \bwy_j, \bx_j + \xi_n \alpha_{ij} \bwy_i \Big\}
\end{eqnarray}
when $J = 0$. These are precisely the subspaces referred to as $W_k^n$ for simplicity in the statement of Theorem~\ref{general_vaaler}.

\begin{lem}\label{Wisotropic}
For each integer $n \geq 1$ and each pair of indices $(i,j) \in \I$, $W^n_{ij}$ is a maximal totally isotropic subspace of $V$.
\end{lem}
\proof
It is clear that $\dim_K W^n_{ij}$ is equal to $\lambda + \W$.  We claim that $F(\bz)=0$ for all $\bz \in W^n_{ij}$.  First assume that $J>0$, then
$$\bz = \bz' + \sum_{k=1,k \neq i}^{\W} a_k \bx_k + a_i \left( \bx_i+\xi_n^2 \alpha_{ij} \bwy_i + \xi_n \bu_j \right),$$
where $\bz' \in V^{\perp}$ and $a_k \in K$ for all $1 \leq k \leq \W$. Recall that $\bz',\bx_1,\dots,\bx_{\W},\bu_j$ are all orthogonal to each other, and $\bwy_i$ is orthogonal to $\bz'$ and all $\bx_k$ with $k \neq i$, as well as $F(\bz') = F(\bx_k) = F(\bwy_k) = 0$ for all $k$. Then we have
\begin{eqnarray*}
F(\bz) & = & a_i^2 F \left( \bx_i+\xi_n^2 \alpha_{ij} \bwy_i + \xi_n \bu_j \right) \\
& = & a_i^2 \left( 2\xi_n^2 \alpha_{ij} F(\bx_i,\bwy_i) + \xi_n^2 F(\bu_i) \right) = 0.
\end{eqnarray*}
The case $J = 0$ can be verified similarly.
\endproof

It is clear from the definition of the $W^n_{ij}$ that if $J > 0$, then for any integers $n, m \geq 1$ and ordered pairs $(i,j), (i',j') \in \I$,
\begin{equation}
\label{dim_W_int_1}
\dim_K \left( W^n_{ij} \cap W^m_{i'j'} \right) = \left\{ \begin{array}{ll}
\lambda+\W-1 & \mbox{if $n \neq m$ and $(i,j) = (i',j')$} \\
\lambda+\max\{0,\W-3 \} & \mbox{if $(i,j) \neq (i',j')$;}
\end{array}
\right.
\end{equation}
whereas, if  $J = 0$, then
\begin{equation}
\label{dim_W_int_2}
\dim_K \left( W^n_{ij} \cap W^m_{i'j'} \right) = \left\{ \begin{array}{ll}
\lambda+\W-2 & \mbox{if $n \neq m$ and $(i,j) = (i',j')$}, \\
\lambda+\W-2 & \mbox{if $(i,j) = (j',i')$}, \\
\lambda+\max\{0,\W-4\} & \mbox{if none of $i,j$ equals any of $i',j'$.} \\
\end{array}
\right.
\end{equation}
In particular, notice that for each ordered triple $(n,i,j)$ we get a different maximal totally isotropic subspace $W^n_{ij}$ of $(V,F)$.

\begin{lem} \label{Wgenerate}
For each integer $n \geq 1$, the family of subspaces $\{W^n_{ij}: (i,j) \in \I\}$ spans $V$ over $K$.
\end{lem}
\proof
For any integer $n \geq 1$, let $W^n$ be the subspace of $V$ spanned by the family $\{W^n_{ij}\}$.  It is clear that $V^{\perp}$ and  $\bx_1,\dots,\bx_{\W}$ are in $W^n$.  If $J > 0$, then all of
$$\xi_n^2 \alpha_{11} \bwy_1 + \xi_n \bu_1,\dots,\xi_n^2 \alpha_{1J} \bwy_1 + \xi_n \bu_J,\dots, \xi_n^2 \alpha_{\W1} \bwy_{\W} + \xi_n \bu_1,\dots,\xi_n^2 \alpha_{\W1} \bwy_{\W} + \xi_n \bu_J$$
are also in $W^n$. If, on the other hand, $J=0$, then $\bwy_1,\dots,\bwy_{\W} \in W^n$.  Therefore it is easy to see that in both cases $W^n$ contains $\lambda+2\W+J=L$ linearly independent vectors, hence $W^n=V$ for each $n \geq 1$.
\endproof

We are now ready to complete the proof of Theorem \ref{general_vaaler} in the case $\W>1$.

\proof[Proof of Theorem \ref{general_vaaler}: $\W > 1$ case]
The first two assertions of the theorem are proved by Lemma \ref{Wisotropic}, Lemma \ref{Wgenerate}, (\ref{dim_W_int_1}), and (\ref{dim_W_int_2}). We are left to show that the height of each $W^n_{ij}$ is bounded above as depicted in the theorem.

First assume $J>0$. By Lemma~\ref{lem_4.7},
\begin{equation}
\label{W_ht_1}
\H(W^n_{ij}) \leq N^{\frac{\delta \W}{2}}\, \H(V^{\perp})\, H(\bx_i+\xi_n^2 \alpha_{ij} \bwy_i + \xi_n \bu_j) \prod_{k=1, k \neq i}^{\W} H(\bx_k),
\end{equation}
where $\delta$ is as in \eqref{delta}. Now, by Lemma~\ref{sum_height},
\begin{equation}
\label{W_ht_2}
H(\bx_i+\xi_n^2 \alpha_{ij} \bwy_i + \xi_n \bu_j) \leq 3^{\delta}\, H(1,\xi_n,\xi_n^2 \alpha_{ij}) \,h(\bx_i) \,h(\bwy_i)\, h(\bu_j).
\end{equation}
Notice that by \eqref{aij},
\begin{eqnarray}
\label{W_ht_2.1}
H(1,\xi_n,\xi_n^2 \alpha_{ij}) & = & H \left( 2F(\bx_i,\bwy_i),2\xi_n F(\bx_i,\bwy_i),\xi_n^2 F(\bu_j) \right) \nonumber \\
& \leq & 2^{\delta}\, h(1,\xi_n^2)\, H \left( F(\bx_i,\bwy_i), F(\bu_j) \right) \nonumber \\
& = & 2^{\delta}\, h(1,\xi_n^2)\, H \left( \sum_{s=1}^N \sum_{t=1}^N f_{st} x_{is} y_{it}, \sum_{s=1}^N \sum_{t=1}^N f_{st} u_{js} u_{jt} \right) \nonumber \\
& \leq & (2 N^2)^{\delta}\, h(1,\xi_n^2)\, h(\bx_i)\, h(\bwy_i) \,h(\bu_j)^2\, H(F).
\end{eqnarray}
Combining \eqref{W_ht_1}, \eqref{W_ht_2}, and \eqref{W_ht_2.1} with \eqref{siegel_for_U}, we obtain:
\begin{eqnarray}
\label{W_ht_3}
\H(W^n_{ij}) & \leq & (6 N^{\frac{\W+4}{2}})^{\delta}\, (C_K(J) \E_K(J)^{1-\delta})^3 \times \nonumber \\
&\ & \times H(F)\, \H(U)^3 \,\H(V^{\perp})\, h(1,\xi_n^2)\, h(\bwy_i)^2\, h(\bx_i)\, \prod_{k=1}^{\W} h(\bx_k).
\end{eqnarray}
Now \eqref{W_ht_3} combined with Lemma \ref{hyper} implies that
\begin{equation}
\label{W_ht_4.1}
\H(W^n_{ij}) \leq \CC'_K(N,\W,J)\, a_K(n)\, H(F)^{\frac{\W+9}{2}} \,\H(U)^3 \,\H(V^{\perp})\, \H(\hyp_i)^7 \,\prod_{k=1}^{\W} \H(\hyp_k),
\end{equation}
where $a_K(n)$ is as in \eqref{akn} if $K$ is a global field, and
\begin{equation}
\label{W_ht_4.2}
\H(W^n_{ij}) \leq \CC'_K(N,\W,J)\, H(F)^{\frac{\W+9}{2}}\, \H(U)^3\, \H(V^{\perp})\, \H(\hyp_i)^{10} \,\prod_{k=1}^{\W} \H(\hyp_k)^2
\end{equation}
if $K=\qbar$. The constant $\CC'_K(N,\W,J)$ in the inequalities above is given by
\begin{equation}
\label{W_ht_4.3}
\CC'_K(N,\W,J) = \left\{ \begin{array}{ll}
9 \cdot 2^{\frac{3\W+19}{2}}\, N^{\frac{\W+12}{2}}\, C_K(J)^3 \, G_K^4\,  B_K(1)^{2\W+6} & \mbox{if $K = $  number field}, \\
16\, q^{\frac{4(\W+3)g(K)}{d}}\, \left( C_K(J) \E_K(J) \right)^3\, G_K^4 & \mbox{if $K = $ function field}, \\
1492992 \cdot 72^{\W+1}\, N^{\frac{\W+12}{2}}\, C_K(J)^3 \, G_K^4 & \mbox{if $K=\qbar$,}
\end{array}
\right.
\end{equation}
and $G_K = \E_K(2)^{1-\delta} A_K(2) C_K(2)$, as in Lemma \ref{hyper}.  For the sake of convenience, let
\begin{equation}
\label{powers_p}
\left\{ \begin{array}{ll}
\p_1(\W) = (\W+5)(\W+1)(\W+2) \left(\frac{3}{2}\right)^\W \\
\p_2(\W) = \frac{(6\W+9) + \p_1(\W)(6\W+5)}{4\W+2}.
\end{array}
\right.
\end{equation}
Now we combine \eqref{W_ht_4.1}, \eqref{W_ht_4.2} with the bounds of Theorem \ref{witt_decomp} to obtain part (iii) of Theorem \ref{general_vaaler}. The constant $\CC_K = \CC_K(N,\W,L,r)$ is defined by
\begin{equation}
\label{W_main_constant_nf}
\CC_K =  \left\{ \begin{array}{ll}
\CC'_K(N,\W,J)\,  \G_K(N,L,\W)^{\W+10}\, a_K(n) \,B_K(r)^r & \mbox{if char$K = 0$,} \\
\CC'_K(N,\W,J)\,  \G_K(N,L,\W)^{\W+10}\, a_K(n) \,q^{\frac{rg(K)}{d}} & \mbox{if char$K > 0$,}
\end{array}
\right.
\end{equation}
when $K$ is a global field, and
\begin{equation}
\label{W_main_constant_qb}
\CC_K = 8\, (2\W+1)^{\frac{3}{2}} \,\CC'_K(N,\W,J) \, \eta(L,r)^{\p_2(\W)}\, 3^{\frac{(10\W+11)(\W+2)\W + 48 \W^4 \p_1(\W)}{2(\W+2)}}\, \W^{\frac{\p_1(\W)}{2}}
\end{equation}
when $K=\qbar$.

Now assume $J=0$, then employing the same kind of estimates as above we obtain:
\begin{eqnarray}
\label{W0_ht_1}
\H(W^n_{ij}) & \leq & N^{\frac{\delta \W}{2}} \,\H(V^{\perp})\, H(\bx_i+\xi_n \bwy_j)\, H(\bx_j+\xi_n \alpha_{ij} \bwy_i)\, \prod_{k=1, k \neq i,j}^{\W} H(\bx_k) \nonumber \\
& \leq & N^{\frac{\delta \W}{2}}\, \H(V^{\perp})\, H(1,\xi_n)\, H(1,\xi_n \alpha_{ij})\, h(\bwy_i)\, h(\bwy_j) \,\prod_{k=1}^{\W} h(\bx_k) \nonumber \\
& \leq & N^{\frac{\delta \W}{2}}\, \H(V^{\perp})\, H(1,\xi_n)^2 \,H\left( F(\bx_i,\bwy_i),F(\bx_j,\bwy_j) \right)\, h(\bwy_i) \,h(\bwy_j) \,\prod_{k=1}^{\W} h(\bx_k) \nonumber \\
& \leq & N^{\frac{\delta \W}{2}}\, \H(V^{\perp})\, H(1,\xi_n)^2\, H(F)\, h(\bwy_i)^2\, h(\bwy_j)^2 \,h(\bx_i)\, h(\bx_j)\,  \prod_{k=1}^{\W} h(\bx_k).
\end{eqnarray}
Then again Lemma \ref{hyper} together with Theorem \ref{witt_decomp} imply part (iii) of Theorem \ref{general_vaaler}. The constant $\CC_K$ in this case is defined by
\begin{equation}
\label{W0_ht_F_const_nf}
\CC_K =  \left\{ \begin{array}{ll}
81 \cdot 2^{\frac{3\W}{2} + 17}\, N^{\frac{\W+16}{2}}\, B_K(1)^{2\W+12} B_K(r)^r \,G_K^8 \,\G_K(N,L,\W)^{\W+14} & \!\!\!\!\mbox{if char$K = 0$} \\
256 \cdot q^{\frac{g(K)(4\W+r+24)}{d}}\, G_K^8 \,\G_K(N,L,\W)^{\W+14} & \!\!\!\!\mbox{if char$K > 0$}
\end{array}
\right.
\end{equation}
when $K$ is a global field, and
\begin{equation}
\label{W0_ht_F_const_qb}
\CC_K = 2^{3\W + 26}\, 3^{\frac{24 \W^4 \q_1(\W)}{\W+2} + \frac{L(L-1)}{2} + 2\W + 16}\, N^{\frac{\W+16}{2}}\ G_K^8 \,\W^{\frac{\q_1(\W)}{2}}\,  \eta(L,r)^{\frac{(6\W+5) \q_1(\W)}{4\W+2}}
\end{equation}
when $K=\qbar$, where $\q_1(\W) = (\W+10)(\W+1)(\W+2) \left(\frac{3}{2}\right)^\W$. This completes the proof of the theorem when $\W>1$.
\bigskip

Now we handle the case $\W=1$.

\proof[Proof of Theorem \ref{general_vaaler}: $\W = 1$ case]
Assume that $\W=1$, then $J \geq 1$ since $r=J+2\W \geq 3$, and we define
\begin{eqnarray}
\label{wnij_3}
& & W^n_1 = \spn_K \left\{ V^{\perp}, \bx_1 \right\},\ W^n_2 = \spn_K \left\{ V^{\perp}, \bwy_1 \right\}, \nonumber \\
& & W^n_{j+2} = \spn_K \left\{ V^{\perp}, \bx_1+\xi_n^2 \alpha_{1j} \bwy_1 + \xi_n \bu_j \right\}
\end{eqnarray}
for each $1 \leq j \leq J$ and $n \in \Nn$. It is clear that these are $(\lambda+1)$-dimensional maximal totally isotropic subspaces of $(V,F)$, and any two of these intersect only in $V^{\perp}$, hence their intersection has dimension $\W+\lambda-1$. Further, for each $n$,  the $J + 2$ subspaces defined in \eqref{wnij_3} span $V$. 

To finish the proof, we need to show that for each $n \in \Nn$ and $1 \leq k \leq J+2$, $\H(W^n_k)$ satisfies the bound of part (iii) of Theorem~\ref{general_vaaler}. By Lemma~\ref{lem_4.7}, for each $1 \leq k \leq J+2$,
\begin{eqnarray}
\label{W1-1}
\H(W^n_k) & \leq & N^{\frac{\delta}{2}}\, \H(V^{\perp})\, \max \left\{ H(\bx_1), H(\bwy_1), H(\bx_i+\xi_n^2 \alpha_{ij} \bwy_i + \xi_n \bu_j) \right\} \nonumber \\
& \leq & (9N)^{\frac{\delta}{2}}\, \H(V^{\perp})\, H(1,\xi_n,\xi_n^2 \alpha_{1j})\, h(\bx_1)\, h(\bwy_1)\, h(\bu_j)
\end{eqnarray}
where $j=k-2$ and the last inequality follows by Lemma~\ref{sum_height}. Combining~\eqref{W1-1} with \eqref{W_ht_2.1}, we obtain:
\begin{eqnarray}
\label{W1-2}
\H(W^n_k) & \leq & (36 N^5)^{\frac{\delta}{2}}\, \H(V^{\perp})\, h(1,\xi_n^2)\, h(\bx_1)^2\, h(\bwy_1)^2\, h(\bu_j)^3\, H(F) \nonumber \\
& \leq & (6 N^{\frac{5}{2}})^{\delta}\, \left( C_K(J) \E_K(J)^{1-\delta} \right)^3\, H(F) \,\H(U)^3\, \H(V^{\perp}) \times\nonumber\\
 & & \times \, \, h(1,\xi_n^2)\, h(\bx_1)^2\, h(\bwy_1)^2,
\end{eqnarray}
where the last inequality follows from~\eqref{siegel_for_U}. Now the statement of the lemma follows from the observation that the bound of~\eqref{W1-2} is precisely the bound of~\eqref{W_ht_3} in this case. This completes the proof of Theorem~\ref{general_vaaler}.
\endproof
\bigskip

{\bf Acknowledgment.} We would like to thank the referee for the helpful suggestions, which  improved the quality of the paper.
\bigskip

\bibliographystyle{plain}  
\bibliography{quad_zero27}
\end{document}